\documentclass[10pt]{amsart}

\usepackage{times}
\usepackage{amssymb}
\usepackage{dbnsymb}
\usepackage{color}
\usepackage{amsfonts,amsmath, amscd}
\usepackage[all]{xy}
\usepackage{graphicx}
\usepackage[margin=24pt,font={small,sf}]{caption}

\theoremstyle{plain}
\newtheorem{theorem}{Theorem}[section]
\newtheorem{lemma}[theorem]{Lemma}

\theoremstyle{definition}
\newtheorem{definition}[theorem]{Definition}

\begin{document}

\author[8]{Fionntan Roukema}
\address{Department of Mathematics, University of Toronto, Toronto, Ontario, M5S 2E4}
\email{fionntanroukema@gmail.com}

\title[10pt]{Goussarov-Polyak-Viro combinatorial formulas for finite type invariants}
\maketitle

\begin{abstract}
Goussarov, Polyak, and Viro proved that finite type invariants of knots are ``finitely multi-local'',
meaning that on a knot diagram, sums of quantities, defined by local information,
determine the value of the knot invariant (\cite{GPV}). The
result implies the existence of Gauss diagram combinatorial formulas for finite
type invariants. This article presents a simplified account of the
original approach. The simplifications provide an easy
generalization to the cases of pure tangles and pure braids. The associated
problem on group algebras is introduced and used to prove the
existence of ``multi-local word formulas'' for finite type invariants of pure braids.
\end{abstract}

\section{Introduction}

\subsection{Statement (Informal)}

\begin{theorem}
A type $n$ invariant $\nu$ can be computed on a knot, represented by a diagram $D$, by studying all subdiagrams of $D$
having up to $n$ crossings.
\end{theorem}

The purpose of section one is to make the informal theorem, formal.

\subsection{Subdiagrams}

For this article, all objects are long, so knots are long knots, pure tangles are long pure tangles etc.

\par
\begin{figure}
\centering
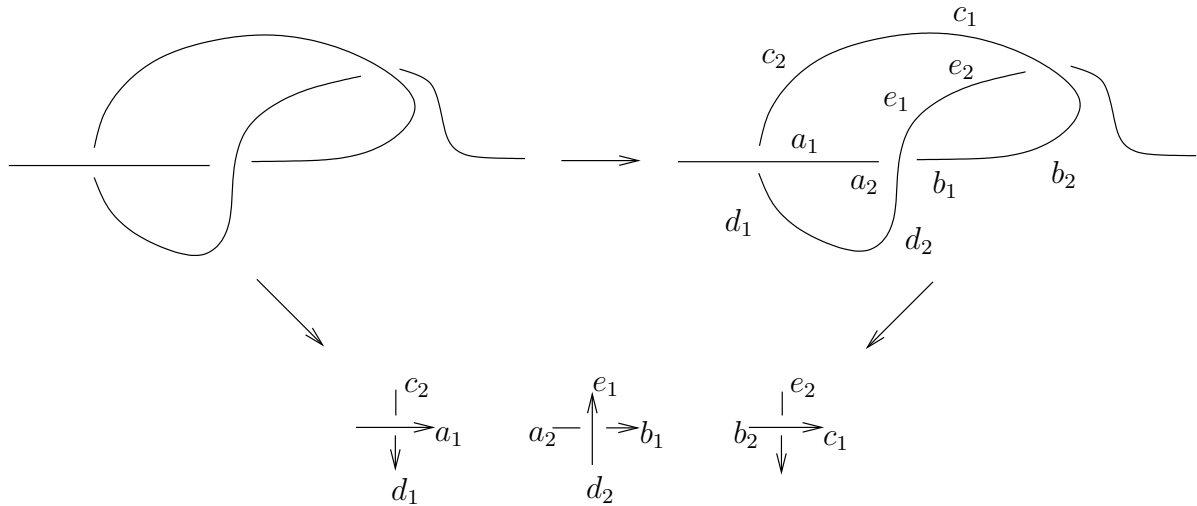
\caption{The result of passing from the long left hand trefoil to
labelled crossing information.} \label{fig1}
\end{figure}

 Our first task is to explain what we mean by ``subdiagrams''. Given a knot diagram, we can restrict our attention to the crossing information.
The connecting arcs between crossings induce a natural labelling on pairs of endpoints of the crossings. See Fig. \ref{fig1}.

\begin{figure}[h]
\centering
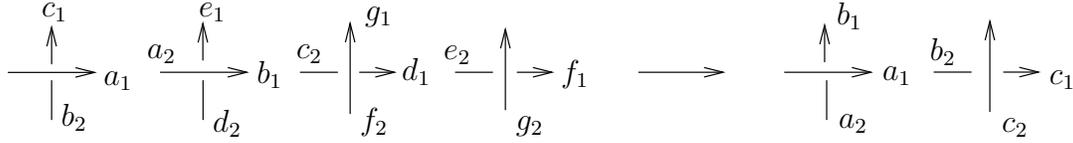
\caption{Passing to the subdiagram containing the second and fourth crossings.}
\label{fig12}
\end{figure}

This motivates the simple idea of thinking of a knot diagram as crossings with
labels on pairs of endpoints explaining how the crossings connect. This synonymous way of thinking about knot diagrams points to a natural
operation, namely, that of considering a subset of crossings with new labelling on pairs of endpoints given by connections
which are now permitted to pass the through deleted crossings.
The result is called a \emph{subdiagram} of the knot diagram. See Fig. \ref{fig12}.

\begin{figure}[h]
\centering
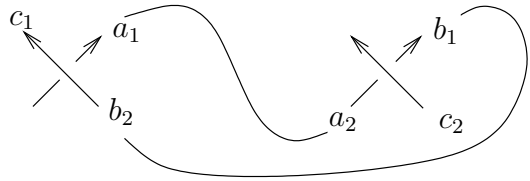
\caption{An example of a connection between the labelled crossings. It is easy to enumerate all possible attempts, all of which
introduce an additional crossing.}
\label{fig2}
\end{figure}

An important point is that a subdiagram of a knot diagram may not give rise to a
real knot diagram. For example, the third connection of Fig. \ref{fig2}
cannot be made without the connecting strand
intersecting some other part of the diagram.

This motivates the definition of a \emph{virtual knot diagram} as
crossings with paired endpoints which may be denoted by labelling, or by a connecting path.

\par We introduce notation for the relevant spaces and maps that
will be used in the paper:

\begin{enumerate}
\item $\mathcal{L}$, the $\mathbb{Z}-$module of formal linear combinations of real long knot diagrams.
\item $\mathcal{VL}$, the $\mathbb{Z}-$module of formal linear combinations of virtual long knot diagrams.
\item The $\mathbb{Z}-$linear map $s:\mathcal{VL} \rightarrow \mathcal{VL}$ which takes a
diagram to the sum of all its subdiagrams, and extended to all of $\mathcal{VL}$ by linearity.
\item $\mathcal{K}$, the $\mathbb{Z}-$module of formal linear combinations of long real knots.
\end{enumerate}

\subsection{Statement (formal)}

We are now in a position to formally state Goussarov's theorem,
and explain its significance.

\begin{theorem}

Let $A$ be an Abelian group. By way of $s$, any $A$-valued type $n$ invariant,
$\nu$, of real knots, considered as a function on $\mathcal{L}$, factors through $\mathcal{VL}$, with the factorization vanishing
on diagrams of $\mathcal{VL}$ with greater than $n$ crossings.

\end{theorem}

In other words there exists a function $\omega : \mathcal{VL} \rightarrow A$, dependant on $\nu$, so that $\omega$
vanishes on virtual knot diagrams with greater than $n$ crossings, and so that the following diagram commutes:

$$
\xymatrix{
{\mathcal{L}} \ar[r]^{\nu} \ar[d]_{s} & A \\
\mathcal{VL} \ar[ur]_{\omega}
}$$

The informal statement is now formal, for the result can be
interpreted as saying ``to compute the value of a knot diagram $D$
under a finite type invariant, it is enough to know $\omega$, and
the subdiagrams of $D$ with less than or equal to $n$ crossings''.
This is the sense in which finite type invariants are finitely
multi-local.
\par Metaphorically, finite type invariants of degree $n$ can be computed with
$n$ fingers;  knowing the values of $\omega$ on all virtual knot diagrams with less than or equal to $n$ crossings,
the $n$ fingered mathematician lets his fingers rest on all combinations of crossings that he can,
and adds the resulting values of the subdiagrams!
\par When $A$ is $\mathbb{Z}$ or $\mathbb{Q}$ the factorization is the existence of Gauss diagram formulas; let $<.,.>$
denote the dirac inner product on $\mathcal{VL}$, then $$\nu(D)  = \omega s(D)= \left \langle \sum \omega(B),s(D) \right \rangle $$
where the sum is over virtual diagrams with less than or equal to $n$ crossings.

\subsection{Content of paper}
Section two provides a proof of the main result, along the lines
of \cite{GPV}, though without invoking Gauss diagrams. The
approach has the effect of allowing some simplifications to the
original proof, and of making generalizations to the context of
pure tangles natural. Section three is devoted to the
exploration of this and other generalizations and related questions.

\section{Proof}

\subsection{Scheme}

Our proof will consist of constructing and explaining the following diagram:

$$
\xymatrix{
\mathcal{K} \ar[dr]^{\nu}\\
{\mathcal{L}}  \ar@<1ex>[d]^{\iota} \ar[r]^{\nu} \ar[u]^{\pi} & A \\
{\mathcal{VL}} \ar@<1ex>[u]^{P} \ar@<1ex>[r]^{s}€€
\ar[ur]^{\overline{\nu}} & {\mathcal{VL}} \ar[u]_{\omega}
\ar@<1ex>[l]^{s^{-1}}
}$$

We will show all triangles in this diagram commutes, and that $\omega$ vanishes on diagrams of order greater than $n$.
This forces $\overline{\nu}$ to be $\nu P$ and $\omega$ to be $\overline{\nu} s^{-1}$.
\par The map $\omega$ is determined by specifying $\overline{\nu}$, and so the correct definition of $P$ becomes the main
issue in the proof. We know the map $P$ must be a well defined means of assigning a real knot diagrams to every virtual knot diagram with real
knot diagrams going to elements in their equivalence class in $\mathcal{K}$. Thus $\pi P \iota= \pi$. The punch line in the proof will be that on certain types of
diagram, the specified $P$ is constant on subdiagrams.
\par The structure of the proof is then:

\begin{enumerate}
\item First, extend $\nu$, considered as a function on
$\mathcal{L}$, to a function
$\overline{\nu}$ defined on $\mathcal{VL}$. This is achieved by
defining $\overline{\nu}$ to be the
composition of $\nu$ with a map $P:\mathcal{VL} \rightarrow
\mathcal{L}$. The definition of $P$
takes some work.
\item Next, $\omega$ is defined to be $\overline{\nu} s^{-1}$.
\item Finally, it remains only to show that $\omega$ vanishes on
elements with greater than $n$
crossings. On ``descending'' diagrams with greater than $n$ crossings
this will be trivial, and we show,
via a map $Q:\mathcal{VL} \rightarrow \mathcal{VL}$, that this is enough.
\end{enumerate}

\subsection{Descendingness}

We start with some definitions that lead to $P$ and thus
$\overline{\nu}$. To do this we will need a notion of
``descendingness'' that will allow us to ``descend'' virtual
crossings in a well defined way. Our notion will require real
crossings to be first encountered as over crossings, with all
double points being ``clumped'' together. This is formalized most
easily in the language of trees. In this setting, $s$ and $s^{-1}$
will still make sense.
\par
\begin{figure}[h]
\centering
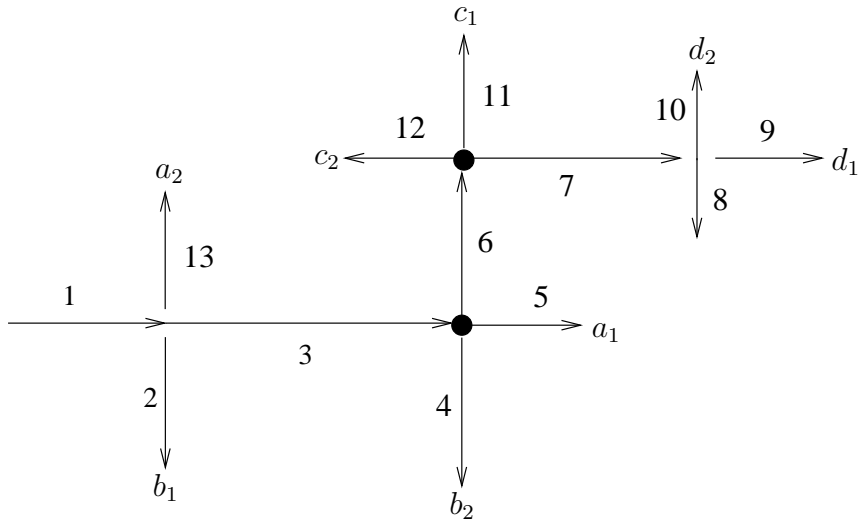
\caption{An example of a tree with information on nodes, arcs and leafs.}
\label{fig13}
\end{figure}

 For us a tree will be directed and carry additional information 
on the nodes, arcs and leafs. The information on nodes will be denoted by double,
over or under crossings. Leafs are grouped in pairs, and arcs are
directed away from the root toward the leafs. The arcs on the tree
will be equipped with ordering given by traversing the tree always
to the ``right'': Starting at the root pass to the right at every
node whenever possible. When no longer possible, turn around and
traverse the tree against the direction until you can turn right.
Eventually the entire tree will be traversed and the arcs are
ordered according to when they are first encountered on the path.
See Fig. \ref{fig13}.

\par
\begin{figure}
\centering
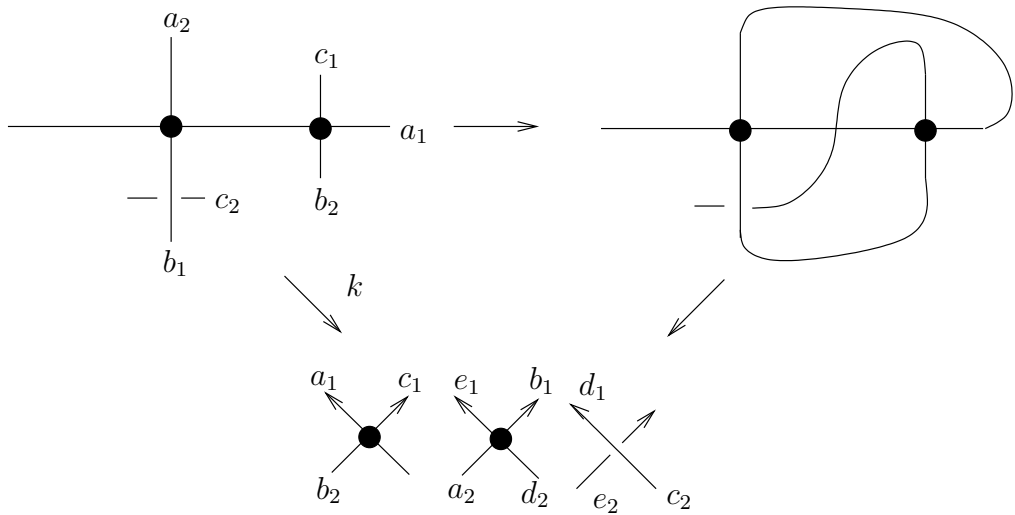
\caption{The map $k$ taking a tree (labelling on arcs and the
directions omitted) to a virtual knot diagram.} 
\label{fig8}
\end{figure}

 Let $k$ be the map from tree diagrams to virtual knot diagrams
given by gluing the ends of the tree together according to the labelling on the leafs.
See Fig. \ref{fig8}.

We call $c$ a \emph{method of finding a tree within a knot diagram} if the following diagram commutes.

$$
\xymatrix{
\parbox{1.7cm}{\centering{virtual knot diagrams}}  \ar@<1ex>[r]^{c} &
\parbox{1.7cm}{\centering{tree diagrams}} \ar@<1ex>[l]^{k}\\
}$$

There is a canonical $c$, and we make it explicit.
\par
\begin{figure}[h]
\centering
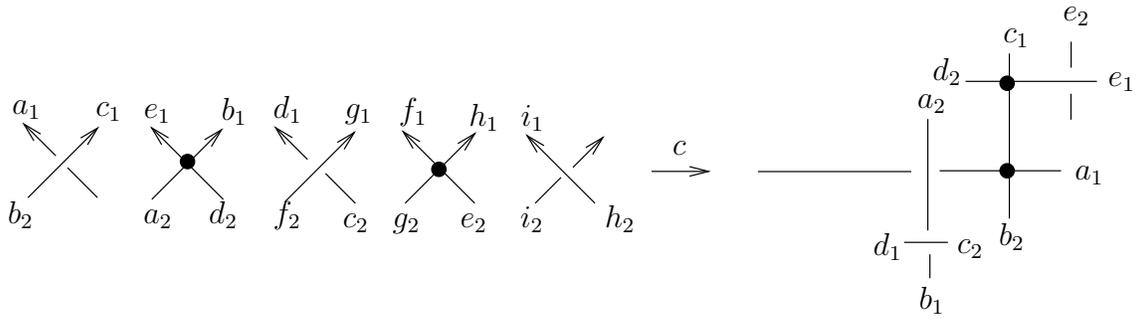
\caption{The map $c$ taking a virtual knot diagram to a tree.}
\label{fig7}
\end{figure}

 Starting with a knot diagram, crossings are glued according
to the labels in ascending order provided homotopy remains trivial
(if homotopy is introduced, the gluing is not performed). The map
$c$ is precisely the map which traverses the knot diagram snipping
the diagram every time it is about to enter a crossing for the
second time. This tree is implicity used in \cite{GPV}. See Fig.
\ref{fig7}.

For our purposes, subtrees having a node of the original tree become a leaf keep the information the node carried.

\begin{definition}
A tree is said to be $\emph{descending}$ if all real crossings
are first encountered as over crossings, and the minimal subtree containing all double points
contains no real crossings.

\begin{figure}[h]
\centering
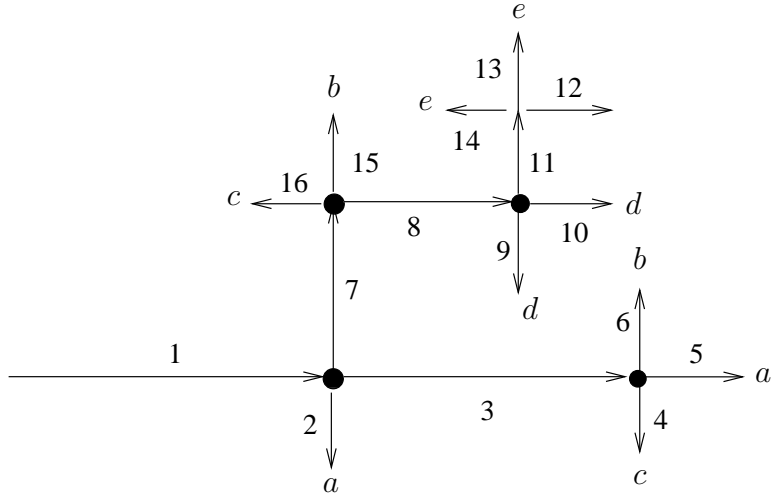
\caption{An example of a descending tree.}
\label{fig3}
\end{figure}

\end{definition}

The maps $c$ and $k$ are extended to formal linear sums of virtual knot diagrams, respectively tree diagrams, by linearity.

\begin{definition}
An element $D$ of $\mathcal{VL}$ is said to be \emph{descending} if $c$ of each element in the summand is descending.
\end{definition}

\subsection{The projection $P$}

We can now define the projection $P$.

\begin{definition}

Given a diagram $D$ from $\mathcal{VL}$, the map $P$ is defined by the following procedure:

\begin{enumerate}

\item Make all real crossings of $c(D)$ descend using the double point relation.
More precisely, if a crossing is first encountered as an under
crossing, we write it as the difference of a double point and an
over crossing using the formal identity $\doublepoint =
\overcrossing - \undercrossing$. This is performed at every
non-descending crossing.

\item For each non-descending diagram, the arcs are directed which gives a notion of left and right
on the arcs (our convention is that ``right'' is taken in the
positive direction). We consider the minimal subtree containing
double points, the arcs of which have induced ordering from the original
diagram. Take the first arc containing real crossings and label
them $r_{1}, \cdots r_{k}$, where $r_{i}$ is to the left of
$r_{j}$ when $i<j$. Isotop the $k^{th}$ crossing through the right
most double point on the arc, then isotop $(k-1)^{st}$ through the
double point etc. The result is a tree with an additional arc free of
real crossings. See Fig. \ref{fig9}.

\begin{figure}[h]
\centering
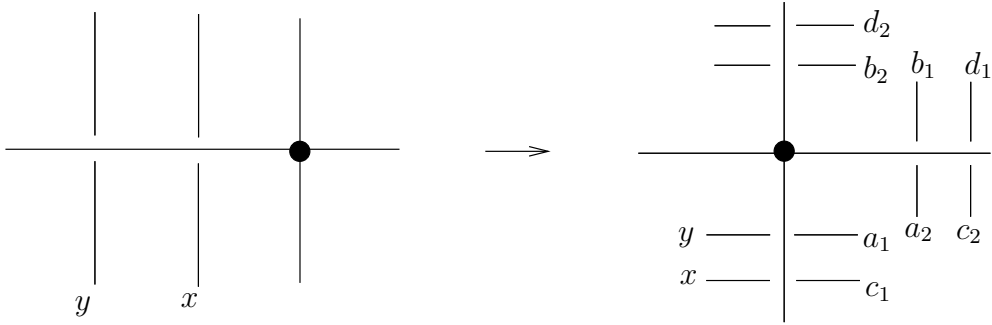
\caption{$x$ and $y$ are markers showing how the tree changes under step two.}
\label{fig9}
\end{figure}

\item Repeat steps one and two until we obtain a sum of descending
diagrams and diagrams with greater or equal to $n$ crossings.
\item Start with the first paired labels and connect them with a path $p_{1}$ off the tree (possible by contractibility of the tree),
then the next labelled pair are connected with a path $p_{2}$ off the tree, and the procedure is followed to yield $t$ paths
$p_{1}, \cdots , p_{t}$. If any two paths $p_{i}$ and $p_{j}$ intersect, with $i<j$, we make the $p_{j}$ descend below $p_{i}$.
This gives an element of $\mathcal{L}$.

\begin{figure}[h]
\centering
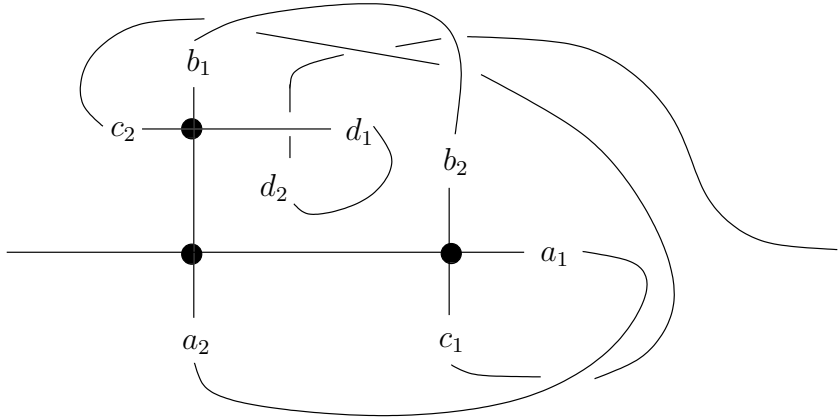
\caption{An example of step four as applied the tree from Fig. \ref{fig3}.}
\label{fig14}
\end{figure}

\item $P$ is extended to all of $\mathcal{VL}$ by linearity.

\end{enumerate}

\end{definition}

Our $P$ serves the same role as the one used in \cite{GPV}, though theirs differs by fixing``bad'' crossings ``one at a time''. Within
our framework it becomes very easy to see that $P$ is well defined and that $P$ reaches step three (see Lemma 2.2).

\subsection{Four Lemmas}

In the above sections $s$ was assumed to be a $\mathbb{Z}-$module isomorphism, and we start by showing it is.

\begin{lemma}
The map $s^{-1}$ exists, and is given by
$$s^{-1}(D):=\sum(-1)^{|D-D'|}D'$$ where the sum is over
subdiagrams $D'$ of $D$, and $|D-D'|$ denotes the difference of the
number of real crossings of $D$ minus the number of real crossings
of $D'$.
\label{lemma1}
\end{lemma}

\begin{proof}
If the diagram has no crossings, then the result is trivial.
Consider the coefficient of an arbitrary element in the summand
of $s^{-1} s(D)$. Any such element $D'$ arising in the sum, is a subdiagram of $D$ and has
coefficient $(1-1)^{|D-D'|}$. This coefficient is zero unless $D'$ equals $D$ in which case the coefficient is one.
\end{proof}

Alternatively, an equally easy symbolic proof:

$$ s^{-1} s( \slashoverback ) = s^{-1} ( \slashoverback + \crossing ) =
\slashoverback - \crossing +\crossing = \slashoverback $$

When $\crossing$ from the symbolic proof is correctly interpreted, Lemma \ref{lemma1} holds.
For example, on a tree $D$, with a specified crossing $C$,
the \emph{subdiagram} without $C$ is given by taking $c$ of the subdiagram of $k(D)$ not containing $k(C)$. Then $\crossing$ is the sum of
subdiagrams not containing the crossing. In particular $s$ and $s^{-1}$ are well defined on trees, braids, and pure tangles.
\par We now show that $P$ is well defined, and that $\overline{\nu}$ extends $\nu$.

\begin{lemma}
$P$ is a well defined map from $\mathcal{VL}$ to $\mathcal{L}$, and $\nu P|_{\mathcal{L}}=\nu$.
\label{lemma2}
\end{lemma}

\begin{proof}
First we check that $P(\doublepoint)=P(\overcrossing - \undercrossing)$, but again this is immediate
by step one in the definition of $P$, for
step one applied to $(\overcrossing - \undercrossing)$ is given by $(\overcrossing - (\overcrossing -\doublepoint))=(\doublepoint)$
meaning $P(\doublepoint)=P(\overcrossing - \undercrossing)$.
\par Each application of step one followed by step two yields a sum of
diagrams with at least the same number of double points and real crossings, each of which have real
crossings off a greater number of arcs of the new tree. Thus, under
repeated application, we eventually find a sum of descending diagrams or diagrams with greater than $n$ crossings.
\par Step four of the operation, involving capping the strands in a descending manner,
is well defined as the subtree of double points determines the
resulting diagram up to Reidemeister moves; any two paths between
the same leafs is off the tree and can be freely isotoped to one
another.
\par Lastly, $P$ does nothing up to Reidemeister moves and so $\nu P|_{\mathcal{L}}=\nu$ completing the proof.
\end{proof}

Next we show that on descending diagrams, $\omega$ has the desired property.

\begin{lemma}
Any descending diagram containing at least one real crossing, or
greater than $n$ crossings is sent to zero by $\omega$.
\label{lemma3}
\end{lemma}

\begin{proof}
The first step is to show that $s$ and $s^{-1}$ preserve double points, but this is clear for
$$s^{-1}(\doublepoint) = s^{-1}(\overcrossing - \undercrossing)=\overcrossing - \crossing - (\undercrossing - \crossing) =
\doublepoint $$  and
$$s(\doublepoint) = s(\overcrossing - \undercrossing)=\overcrossing + \crossing - (\undercrossing
+ \crossing) =\doublepoint $$
Now suppose $D$ has greater or equal to $n$ crossings. If greater than $n$
of the crossings are double points then we are done by finite typeness of $\nu$. Otherwise $D$ must
contain at least one real crossing. So
$$\omega(D)=\sum(-1)^{|D-D'|}\overline{\nu}(D')$$
But $D$ descending means all subdiagrams have the same minimal double point tree
and consequently that $\overline{\nu}(D')=\overline{\nu}(D)$ for every subdiagram $D'$.
Whence $$\omega(D) =(\sum(-1)^{|D-D'|})\overline{\nu}(D)=(1-1)^{k}$$
where $k$ is the number of real crossings of $D$ which we assumed to be greater than zero implying $\omega(D)=0$
as required.
\end{proof}

In the course of the proof we showed that $s$ preserved double points and thus that $s$ of a descending diagram is descending.
This will be used again.

\begin{lemma}
Any diagram $D$ with greater or equal to $n$ crossings is sent to zero
by $\omega$.
\label{lemma4}
\end{lemma}

\begin{proof}
The case where $D$ has no real crossings is trivial,
so we assume that $D$ contains at least one real crossing. Set
$$Q=s\iota P s^{-1}$$ As $P s^{-1}(D)$ is descending we see $Q(D)$ is descending.
Now $$\overline{\nu} s^{-1} Q:= \nu P  s^{-1} s \iota P  s^{-1}=\nu  P
\iota P s^{-1}= \nu P  s^{-1}=\overline{\nu} s^{-1}$$
We have already seen $P$, $s$, and $s^{-1}$ all maintain at least the same number
of double points, and that $Q(D)$ is descending, meaning that we're done if each element of
$Q(D)$ contains at least one real crossing (see Lemma \ref{lemma3}).
\par We set $p$ as a single application of steps one and two from the
definition of $P$, and $d$ as step four, then $P$ is given by
$dp^{m}c$ for some integer $m$. Now $$s \iota P s^{-1}(D)=s \iota d
p^{k} s^{-1}= s \iota d p s^{-1} s p s^{-1} \cdots s p cs^{-1}$$
Adopting the notation of Fig. \ref{fig6}, we considering a single conjugate of $p$. The only subdiagram coming
from $p(D_{2})$ with no real crossings cancels with the only
subdiagram from $p(D_{1})$ with no real crossings, for the
diagrams are equivalent up to virtual Reidemeister three
(alternatively, virtual crossings are off the tree), and have opposite signs. Thus a single
conjugate of $p$ keeps at least one real crossing.

\begin{figure}[h]
\centering
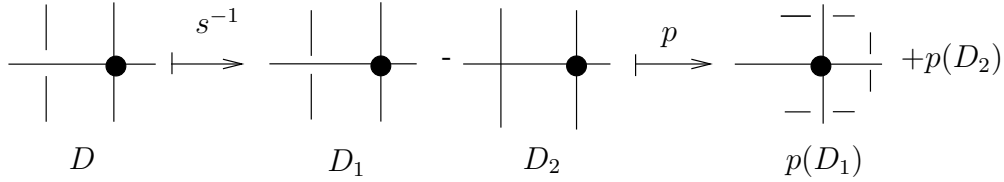
\caption{$D_{1}$ represents the sum of subdiagrams which contain the real crossing depicted, $D_{2}$
represents with the sum of subdiagrams which don't contain the real crossing.}
\label{fig6}
\end{figure}

If we replace $p$ with $pc$ in Fig. \ref{fig6}, we reason in the same way to see $spcs^{-1}$ keeps one real crossing.
\par Write $s \iota d p s^{-1}(D)=d s p s^{-1}(D)+B$ where each element of $B$ is descending with at least one real crossing
coming from the descending caps. We see all elements of $Q(D)$ to be descending and to have at least one real crossing.
By Lemma  \ref{lemma3} the sum must go to zero under $\omega$ as claimed.

\end{proof}

In \cite{GPV} the operator $Q$ is close to $d s p^{m} c s^{-1}$, and the argument requires more work.

\subsection{Putting the pieces together}

We have proved Goussarov's theorem for we have shown
$\omega(D) = \omega Q (D)$ (Lemma \ref{lemma4}) and that $Q(D)$ is killed by $\omega$ when $D$ has greater than
$n$ crossings (Lemma \ref{lemma4}), and

$$
\xymatrix{
{\mathcal{L}}  \ar@<1ex>[d]^{\iota} \ar[r]^{\nu}
        & A \\
{\mathcal{VL}} \ar@<1ex>[u]^{P} \ar@<1ex>[r]^{s}
\ar[ur]^{\overline{\nu}} & {\mathcal{VL}} \ar[u]_{\omega}
\ar@<1ex>[l]^{s^{-1}}
}$$

commutes (Lemmas \ref{lemma1}, \ref{lemma2}) which was what we wanted.

\section{Generalizations, and the question in the framework of algebras}

\subsection{Pure tangles}

\par For pure tangles it will be possible to pass to a canonical tree with an
induced ordering, and from a tree back to a pure tangle,
namely, that we have a corresponding $c$ and $k$. By way of the new
$c$ and $k$ we inherit descendingness, $\omega$, $P$,
and $Q$. We may then read the earlier lemmas, theorems, definitions and proofs
with the words ``knots'' read as ``pure tangles'', $\mathcal{L}$ replaced with $\mathcal{PT}$
(the $\mathbb{Z}$-module generated by pure tangle diagrams), and
$\mathcal{VL}$ replaced with $\mathcal{VPT}$ (the $\mathbb{Z}-$module generated by virtual
pure tangle diagrams) to obtain a proof of Goussarov's theorem in the new context of pure tangles.
\par
\begin{figure}[h]
\centering
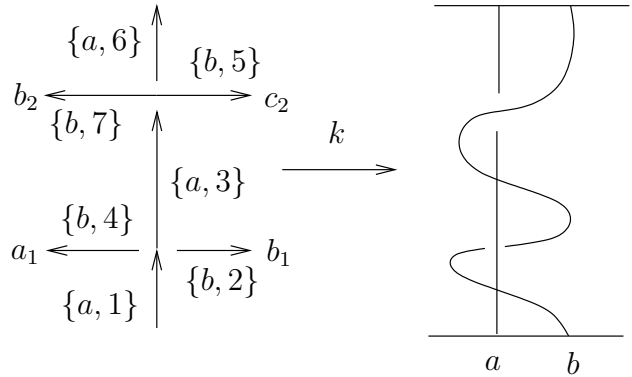
\caption{The effect of taking $k$ of a tree.}
\label{fig11}
\end{figure}

 A \emph{tree} will be as previously defined, but with additional information on the arcs indicating the strand number
(for example, a strand color).  As before, the definition of $k$
is the gluing of the strands according to the labelling on pairs
of leafs. Strand color indicates how to glue to the shield of the
tangle. See Fig. \ref{fig11}.

\begin{figure}[h]
\centering
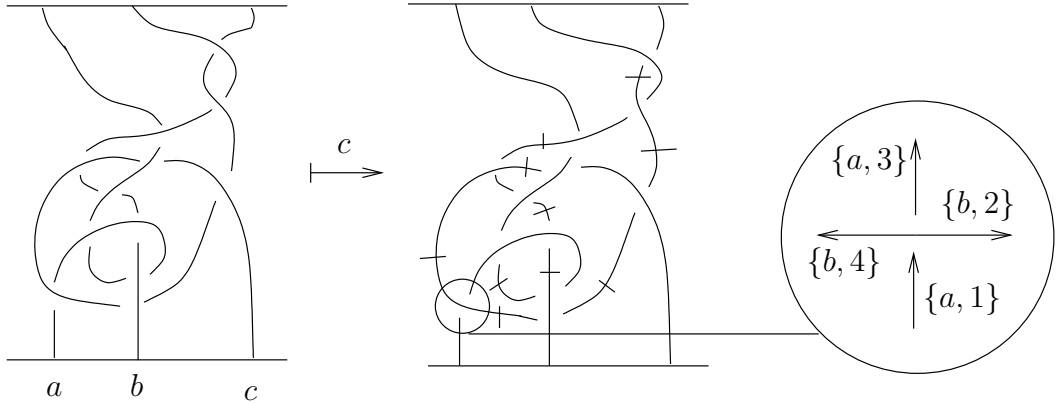
\caption{The enlarged crossing gives an example of the labelling on arcs. The straight line marks on the right
diagram denote snips which should be interpreted as breaks within the diagram to give a tree.}
\label{fig10}
\end{figure}

We exhibit a $c$ to complete the argument. For a diagram $D$, the
map $c$ follows the first strand marking crossings as encountered,
and cutting just before entering a crossing for the second time.
Next, the second strand is followed and cut just before meeting
what has already been traversed. The procedure is followed for the
remaining strands until we arrive at a tree. The tree is directed
as before, the node information is inherited from crossing type,
arcs and leafs are labelled as before but with the extra piece of
information on arcs representing strand number. See Fig.
\ref{fig10} and Fig. \ref{fig11}.

So we have
$$
\xymatrix{
{\mathcal{PT}} \ar[r]^{\nu} \ar[d]_{s} & A \\
\mathcal{VPT} \ar[ur]_{\omega}
}$$

With $\omega$ vanishing on virtual pure tangles of degree less
than or equal to $n$, which is precisely Goussarov's finite
multi-localness theorem for pure tangles.

\subsection{Braids}

In \cite{DBN} it was shown that any finite type invariant of pure
braids can be extended to a finite type invariant, of the same
type, of pure tangles\footnote{The paper uses the terminology
``string link'' for pure tangle.}. Thus we can factor a finite type
invariant of braids through the formal algebra of virtual pure
tangles with the factorization vanishing on virtual pure tangles
with greater than $n$ crossings.

$$
\xymatrix{
PB \ar[r]^{\nu} \ar[d]_{s} & A \\
\mathcal{VPT} \ar[ur]_{\omega}
}$$

This is a Goussarov finite multi-localness type theorem for pure
braids and points to the natural question of whether or not finite
type invariants can be factored through some algebra, via $s$,
with the factorization vanishing on elements of ``order n''? It is
natural to ask this question of group algebras, for which we have
a very simple answer.

\subsection{Finite type invariants on group algebras}

\par
\begin{figure}[h]
\centering
\includegraphics{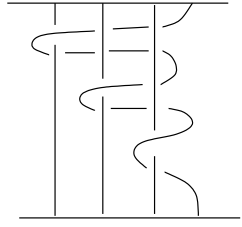}
\caption{The word $abc$ in $F_{3}$ viewed as an element of $PB_{4}$.}
\label{fig4}
\end{figure}

 Let $G$ denote a group, and $RG$ a group algebra generated by elements of $G$. As in the
case of knots, we will need to work with elements in the equivalence class of elements of $RG$. Let $\widetilde{RG}$ denote
the vector space of formal linear combinations of words. Note that $\widetilde{RG}$ includes non-reduced words, and
so though $gg^{-1}$ and $1$ are considered as the same element of $RG$,
they are different elements of $\widetilde{RG}$. We have a natural notion of \emph{finite typeness} from the theory of braids
when $G$ is a free group; an arbitrary word can be thought of as a pure braid with $k+1$ strands, $k$ of
which are straight. Figure \ref{fig4} shows how to think of the word
$abc$ in $F_{3}$ as an element of $PB_{4}$.

\begin{figure}[h]
\centering
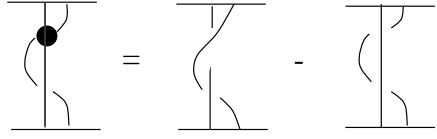
\caption{A double point in $PB_{2}$.}
\label{fig5}
\end{figure}

From which, one sees the appropriate criteria for a function
on $F_{k}$ to be of \emph{ finite type}; we say an invariant
is of type $n$ if it vanishes on words containing greater than
$n$ products of the form $(g-1)$, see Fig. \ref{fig5}.

A function $\nu$ vanishing on words with greater than $n$
subfactors of the form $(g-1)$ is equivalent to $\nu$ vanishing on
words with greater than $n$ subfactors of the form $(g-h)$ because
$(gh^{-1}-1)h=(a-h)$ and $\nu$ is invariant. Note this definition
corresponds to the old notion of finite typeness in the case of
pure braids. A Goussarov finite multi-localness type theorem
holds.

\begin{theorem}
Type $n$ invariants on $\widetilde{RG}$ factor through
$\widetilde{RG}$ with the factorization vanishing on words of
$\widetilde{RG}$ of length greater than $n$.
\end{theorem}

\begin{proof}
The maps $s$ and $s^{-1}$ are defined from $\widetilde{RG}$ to $\widetilde{RG}$ as in sections 1.3 and 2.3 but with ``subdiagrams'' replaced
with ``subword''. Then tautologically we have a commutative diagram:

$$
\xymatrix{
\widetilde{RG} \ar[r]^{\nu} \ar[d]_{s} & A \\
\widetilde{RG} \ar[ur]_{\omega}
}$$

Now $$s^{-1}(g_{j_{1}}g_{j_{2}}\cdots g_{j_{m}})=(g_{j_{1}}-1)(g_{j_{2}}-1)\cdots (g_{j_{m}}-1)$$ and
$$\nu (g_{j_{1}}g_{j_{2}}\cdots g_{j_{m}})=\nu s^{-1} s(g_{j_{1}}g_{j_{2}}\cdots g_{j_{m}})$$ meaning $\nu s^{-1}$
is zero on words with greater than $n$ letters. Thus $\nu$, viewed as a
function on $\widetilde{RG}$, is determined by its values on subwords of length
less than or equal to $n$.
\end{proof}

\subsection{Pure braids again}

On pure braids, the notions of finite typeness on diagrams and
group algebras coincide for the definition on group algebras was
made to emulate the notion on pure braids. Considering $PB_{n}$ as a
finitely presented group, get the following commutative diagram with $\omega$
vanishing on words of length greater than $n$:

$$
\xymatrix{
\widetilde{PB_{n}}  \ar[d]^{s} \ar[r]^{\nu} & A \\
\widetilde{PB_{n}} \ar[ur]_{\omega := \nu s^{-1}}
}$$

This is a finite multi-localness result, giving combinatorial
formulas in terms of subwords.

\subsection{Where next?}

The theorem yields nice combinatorial formulas, and the proof is
entirely constructive. This leads us to first wonder what steps
would be involved in the implementation of a program that computes
the Gauss diagram formula of a finite type invariant? In
principle, this involves translating the steps of the proof into
programming syntax. It is easy to see that the program must only
be capable of computing the tree of a knot diagram (so we need a
specific $c$), and $\omega$ of an arbitrary basis element of
$\mathcal{VL}$ containing less than or equal to $n$ crossings.
\par Understanding $\omega$ involves only knowing how $\overline{\nu}$ and $s^{-1}$ work. For the latter function,
we need to be able to tell a program how to find all subdiagrams of a knot. For the former, we must tell the
program how to take $P$ of an element. In other words we must be able to tell a program how to add
diagrams to make all real crossings descend, and how move a real
crossings to the right.
\par The next questions are when do combinatorial formulas involving a
sum of virtual knots give rise to an invariant of long knots, and
is our extension $\overline{\nu}$ is an invariant (of finite type)
of virtual knots? As observed in \cite{GPV}, the first part of the
first question has an easy answer, namely, we need the formula to
satisfy $s$ of the real Reidemeister moves. The second part of the
question is conjectured to be true in \cite{GPV}, and by the same
reasoning we must check only that the Gauss diagram formula that
we get for $\nu$ satisfies $s$ of the real and virtual
Reidemeister moves. Our approach shows trivially that Reidemeister
moves containing only virtual crossings do not effect
$\overline{\nu}$ (the virtual crossings are considered to be off
the tree, and the subtree of double points determines the value of
$\overline{\nu}$). The implementation of a computer program requires finding $\overline{\nu}$,
which could be used on equivalent virtual knot diagrams
with a view to disproving or supporting the conjecture.
\par Lastly, the simpleness of the result in the context of algebras and its application to pure braids is
pleasing and begs the question as to whether this approach may be extended to other topological
objects?

\subsection{Acknowledgements}
This paper grew out of a graduate reading course with Professor
Dror Bar-Natan at the University of Toronto. Without his help and
guidance, both mathematically and financially, this paper would
not have been written. I would also like to thank Karene Chu, a
Ph.D. student at the University of Toronto, who contributed to the
project. I thank them both sincerely.

\end{document}